\documentclass[final]{IEEEtran}
\usepackage{xcolor,soul,framed} 
\usepackage[pdftex]{graphicx}
\usepackage[switch]{lineno}
\usepackage[cmex10]{amsmath}
\usepackage{array}
\usepackage{mdwmath}
\usepackage{mdwtab}
\usepackage{eqparbox}
\usepackage{url}
\usepackage{tabularx}
\usepackage{lipsum,booktabs}
\usepackage{float}
\usepackage{threeparttable}
\usepackage{multirow}
\usepackage{makecell}
\usepackage{tikz}  
\usepackage[caption=false]{subfig}
\usepackage[T1]{fontenc}
\usepackage[noadjust]{cite}
\usepackage{algorithm}
\usepackage{algpseudocode}
\usepackage{amsmath}
\usepackage{hyperref}
\usepackage{etoolbox}      
\makeatletter
\newcommand{\algmarginnote}[1]{%
  \vskip\ALG@thistlm
  \hrule
  \smallskip
  \textit{\footnotesize #1}
}
\makeatother

\newcolumntype{L}[1]{>{\raggedright\arraybackslash}m{#1}}
\newcolumntype{R}[1]{>{\raggedleft\arraybackslash}m{#1}}
\newcolumntype{C}[1]{>{\centering\arraybackslash}m{#1}}

\begin{document}
	
\title{A Global-Local Optimization Approach for Asynchronous SAR ADC Design}
\author{
Yijia Hao, \emph{Student Member, IEEE}, Ken Li, \emph{Student Member, IEEE}, Miguel Gandara, \emph{Member, IEEE}, Shaolan Li, \emph{Senior Member, IEEE}, Bo Liu, \emph{Senior Member, IEEE}
 \thanks{Y. Hao, S. Cochran, and B. Liu are with School of Engineering, University of Glasgow, Scotland. (e-mails: 2357677h@student.gla.ac.uk, \{Sandy.Cochran, Bo.Liu\}@glasgow.ac.uk)}
 \thanks{K. Li and S. Li are with Georgia Institute of Technology, USA. (emails: kli416@gatech.edu, \ shaolan.li@ece.gatech.edu)}
 \thanks{M. Gandara is with Mediatek Inc. USA. (e-mail: Miguel.Gandara@gmail.com)}
 }
\maketitle

\begin{abstract}
This paper presents a system-level optimization framework for automated asynchronous SAR ADC design, addressing the limitations of block-level methods in terms of suboptimal performance and manual effort. The proposed approach integrates a fast global optimizer with a multi-fidelity local optimizer to efficiently handle high-dimensionality and expensive simulation cost. Experimental results from 12 design cases, covering 7- and 12-bit resolutions and a frequency range of 100 kHz to 250 MHz, demonstrate highly competitive performance compared with prior works.
\end{abstract}

\begin{IEEEkeywords}
successive-approximation-register (SAR), analog-to-digital converters (ADC), design automation, analog IC sizing, optimization
\end{IEEEkeywords}

\section{Introduction\label{sec:introduction}}

As a key interface between the analog and digital domains, successive-approximation-register (SAR) analog-to-digital converters (ADCs) are widely adopted in modern electronic systems for their energy efficiency, moderate-to-high resolution, and suitability for low-power applications. Despite these advantages, SAR ADC design remains a complex task due to the intricate interactions among its core building blocks, including the comparator, digital-to-analog converter (DAC), sample-and-hold (S/H) circuit, and digital control logic. Each of these components introduces specific design constraints that must be carefully balanced to meet overall performance goals. Traditional manual approaches depend heavily on designer experience and iterative tuning, making the process labor-intensive and time-consuming. This has motivated a growing research for automated optimization techniques to streamline SAR ADC design.

Several automated or semi-automated design methods have been proposed to address the need for more efficient SAR ADC design. One approach replaces analog building blocks with synthesizable digital standard cells, reducing manual efforts and improving design portability across technology nodes \cite{seoReusableCodeBasedSAR2018}. However, this restricts the design flexibility as digital standard cells may not always meet stringent ADC performance requirements. Another notable method is hybrid automation, exemplified by \cite{dingHybridDesignAutomation2018, liuOpenSAROpenSource}, which combine techniques such as optimization algorithm, lookup table, and library-based selection for specific blocks. In these methods, designers initially allocate specifications to each individual block, then optimize them independently without fully considering inter-block interactions. As a result, overall performance heavily depends on initial specification allocation. Another systematic method determines block-level component sizing in cycles, which explicitly considering inter-block effects \cite{huangSystematicDesignMethodology2016}. Although system-level trade-offs are considered, the approach still requires significant manual effort, especially for sizing individual blocks. As a result, current automated methods face two major limitations: they remain restricted to block-level design, leading to suboptimal system performance, and they demand significant manual effort when applied to different specifications or technology nodes.

\begin{figure}[t]
\centering
\includegraphics[width=3.2in]{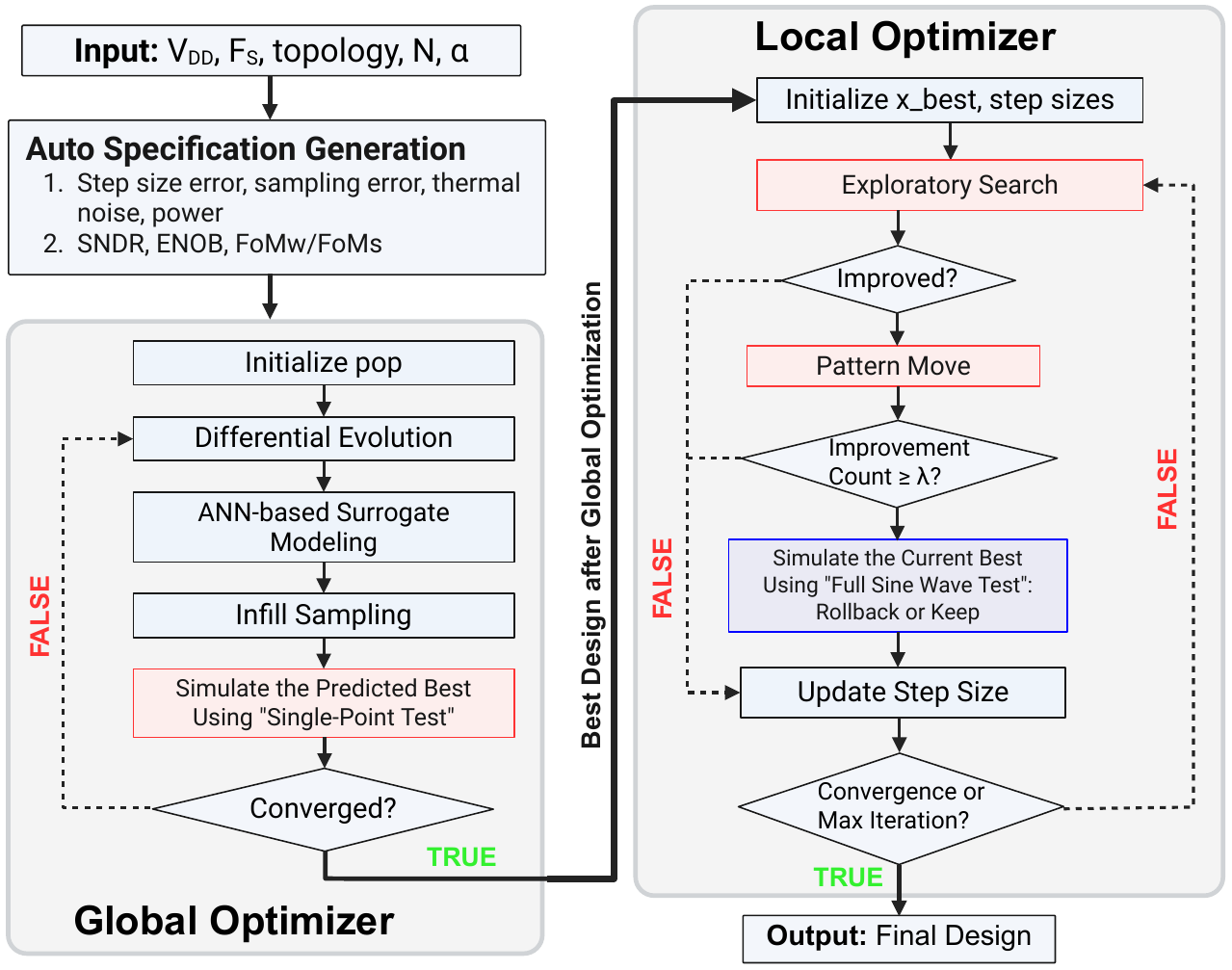}
\caption{The flow diagram of the proposed global-local sizing approach. The red blocks are based on the single-point test, while the blue block represents the full sine wave test.}
\label{flow}
\end{figure}

To address the limitations of prior block-level methods, a system-level design approach is needed. While global optimization offers a holistic solution, its practical use is limited by the long simulation time and high-dimensional design complexity \cite{lyuStudyExploringExploiting2024}. In response to these challenges, this work proposes an efficient global-local optimization framework. The framework integrates a computationally inexpensive global search algorithm (executed within four hours) for broad design space exploration, with a local optimizer that applies parallel, multi-fidelity pattern search to unconverged variables (completed within three hours). During global exploration, performance constraints such as noise, distortion, and power are automatically derived from top-level requirements and evaluated using low-cost single points tests, significantly reducing simulation overhead. Once most of the design variables have converged, the local optimizer is applied to the remaining parameters to accelerate convergence. To ensure system-level performance accuracy, the local optimizer targets the remaining variables using blended objective function and full sine wave test, accelerated by parallel simulation. By combining the two stages, the proposed approach ensures both computational efficiency and high-performance design. Different from prior works which focus on optimizing individual building blocks, this work treats the SAR ADC as an integrated system during sizing, which enables more effective exploration of design trade-offs. 



The remainder of this paper is organized as follows. Section II introduces detailed SAR implementation and design considerations. Section III describes the proposed global-local optimization method. Section IV shows the experiment results. Concluding remarks are presented in Section IV.

\section{Architecture And Design Considerations}

\begin{figure}[t]
\centering
\includegraphics[width=3in]{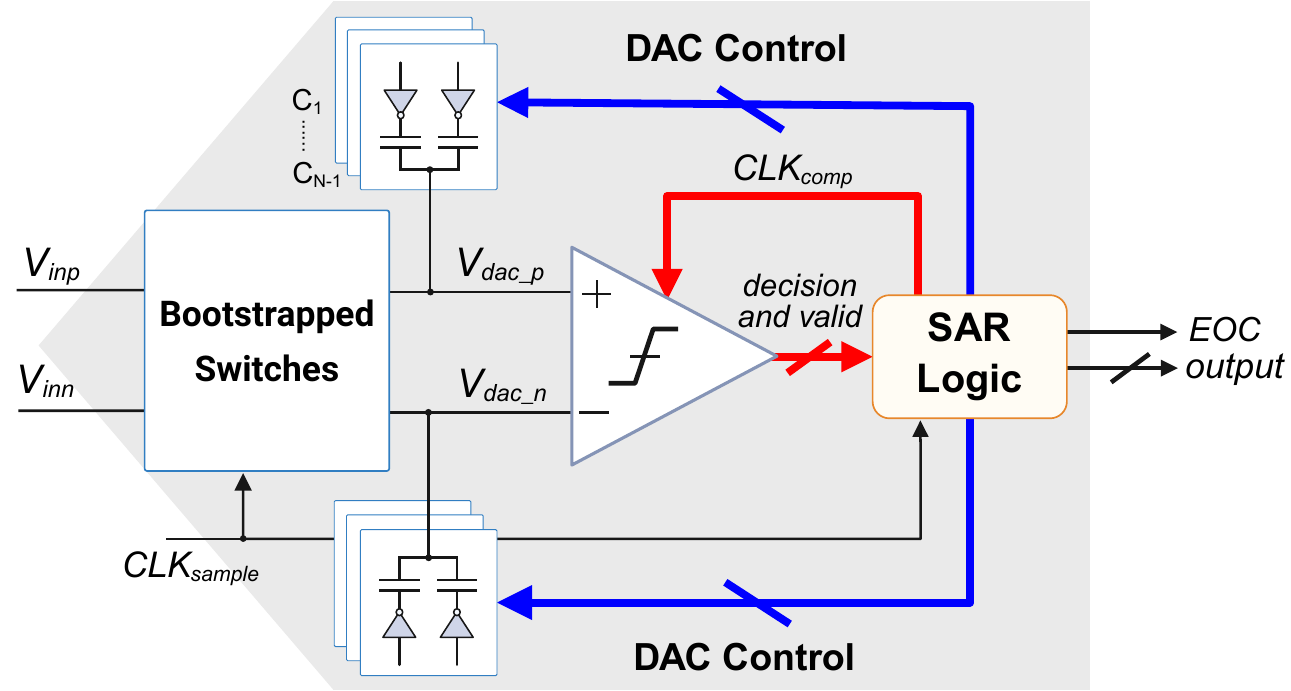}
\caption{The architecture of an N-bit asynchronous SAR ADC.}
\label{architecture}
\end{figure}

Fig. \ref{architecture} shows the selected N-bit differential SAR ADC. In this paper, the topologies of the S\&H, comparator, and DAC are fixed to a bootstrapped switch, a dynamic comparator, and a top-plate sampling and fully binary-weighted C-DAC with a Vcm-based switching scheme \cite{zhu10bit100MSReferenceFree2010a}. These topologies are selected for their wide applicability and ability to support a broad range of speed and resolution requirements. The SAR logic includes a shift register for generating asynchronous clock signals and registers for bit decisions, implemented using dynamic D flip-flops (DFF) for their high energy efficiency \cite{shaikhHighSpeedLow2018}. 

Accuracy, speed, and power are the major trade-offs in SAR ADC design. Two major distortion sources are considered in the optimization process: (1) incomplete charge transfer during sampling due to limited time, and (2) timing misalignment between DAC settling and comparator decisions in asynchronous SAR ADCs, which can lead to settling errors or even conversion failure. Capacitor mismatch is another source but is excluded as it is typically calibrated \cite{ wangSimpleHistogramBasedCapacitor2020, duCapacitorMismatchCalibration2023}. Thermal noise from active devices in the S/H and comparator stages limits SNDR in high-resolution designs and must be carefully managed. Timing is equally critical, as asynchronous SAR ADCs require precise inter-block coordination in the absence of a global clock. All of these considerations are inherently captured in the proposed system-level approach, without the need for manual specification allocation across individual blocks.

\section{Global-Local Optimization Framework}

\subsection{Overview}
The advancement of optimization-based methods and parallel computing has motivated the use of this approach. The SAR ADC sizing is formulated as a constrained optimization problem, mathematically stated as (1).

\begin{equation}
\begin{aligned}
\min_{\mathbf{x}} \quad & f(\mathbf{x}) \\
\text{s.t.} \quad & g_m(\mathbf{x}) \geq G_m \quad \text{or} \quad g_m(\mathbf{x}) \leq G_m \\
& \forall m = 1, 2, \dots, M.
\end{aligned}
\end{equation}

Here, $\boldsymbol{x}$ is the vector of design variables (e.g., capacitor size, transistor size, sampling time), $f(\boldsymbol{x})$ denotes the design objective (e.g., power consumption), and $g_m(\boldsymbol{x})$ represents the $m$-th performance constraint (e.g., sampling error). $G_m$ is the specification target for the $m$-th constraint.

\subsection{Specification Derivation}
Conventional methods decouple sizing from system evaluation. The proposed approach integrates system-level performance from the outset. Table \ref{spec} defines two sets of specifications for each design phase. During the global optimization phase, the specifications considered include step size ratio error (SSRE), sampling error, thermal noise, and power consumption. The parameter $\alpha$ serves as a user-defined scaling factor to flexibly tighten or relax SNDR, enabling trade-off control between performances. $\alpha$ is set to 1 by default. 

In a binary-weighted SAR ADCs, the actual step size for bit $i$ can deviate from ideal due to capacitor mismatch, incomplete settling, parasitics, and other nonidealities. The actual step can be modeled as:
\begin{equation}
\text{step}_i = A_i (1 + \delta_i), \quad \text{with } \frac{A_i}{A_{i+1}} = 2,
\end{equation}
where $\delta_i$ represents the relative error of the $i$th DAC step, capturing both static mismatch and dynamic errors such as incomplete settling and parasitic-induced distortion. Then, the SSRE of two succeeding bits can be defined:
\begin{equation}
\text{SSRE}_i = \left| \frac{\text{step}_i}{\text{step}_{i+1}} - 2 \right| = \left| \frac{2(1 + \delta_i)}{1 + \delta_{i+1}} - 2 \right|.
\end{equation}
For small relative errors (i.e., $\delta_i, \delta_{i+1} \ll 1$), a first-order approximation gives:
\begin{equation}
\text{SSRE}_i \approx 2 |\delta_i - \delta_{i+1}|.
\end{equation}

Assume that the dynamic errors are uncorrelated across bits. This assumption holds approximately when bit-level settling is sufficiently fast to prevent significant error propagation across cycles. To ensure that the dynamic errors do not exceed the quantization noise power, the total error voltage power is constrained as:
\begin{equation}
\sum_{i=1}^{N} V_{\varepsilon_i}^2 < \left( \frac{\Delta}{\sqrt{12}} \right)^2,
\end{equation}
where the error voltage for bit $i$ is:
\begin{equation}
V_{\varepsilon_i} = 2^{N - i} \cdot \Delta \cdot \delta_i.
\end{equation}
If the same error budget is applied for each bit:
\begin{equation}
V_{\varepsilon_i}^2 = \frac{\Delta^2}{12N} \quad \Rightarrow \quad \delta_i = \frac{1}{2^{N-i} \cdot \sqrt{12N}}
\end{equation}
Substituting (7) into the SSRE approximation gives:
\begin{align}
\text{SSRE}_i &\approx 2 |\delta_i - \delta_{i+1}| \notag\\ 
&= 2 \cdot \left| \frac{1}{2^{N - i} \sqrt{12N}} - \frac{1}{2^{N - i - 1} \sqrt{12N}} \right| \notag\\
&=  \frac{1}{2^{N - i - 1} \cdot \sqrt{12N}}
\end{align}


The constraints for sampling error and thermal noise can also be derived under the assumption that their power is equal to the quantization noise. For example, for sampling error, 
\begin{equation}
\text{Sampling Error} = \left| V_{\text{ideal}} - V_{\text{sampled}} \right| < \frac{\Delta}{\sqrt{12}} = \frac{V_{\text{DD}}}{2^N \cdot \sqrt{12}}.
\end{equation}

The upper bound for the SNDR can then be calculated, considering quantization noise, thermal noise together with the SSE and sampling error: 
\begin{align}
P_{\text{tot, max}} &= P_{\text{q}} + P_{\mathrm{thermal,rms}}+ P_{\mathrm{sse}}+ P_{\mathrm{sample}} \notag\\
& = \frac{\text{LSB}^2}{12} + \frac{\text{LSB}^2}{12} + \frac{\text{LSB}^2}{12} + \frac{\text{LSB}^2}{12} = \frac{\text{LSB}^2}{3}
\end{align}

The maximum SNDR is therefore: 
\begin{equation}
\mathrm{SNDR} < 10 \log_{10} \left( \frac{P_{\mathrm{signal}}}{P_{\mathrm{tot, max}}} \right) = 6.02N - 4.25\ \text{dB}.
\end{equation}

\begin{table}[t]
\caption{Summary of Specifications Used in Optimization}
\label{spec}
\renewcommand{\arraystretch}{1.6}
\centering
\begin{threeparttable}
\setlength{\tabcolsep}{1pt}{
\begin{tabular}{p{1.8cm} p{2.5cm} p{4.1cm}}
\hline
\textbf{Level} & \textbf{Performance} & \textbf{Specification} \\
\hline
\rule{0pt}{3.6ex}
\multirow{4}{*}{\makecell[c]{Coarse \\ Evaluation}} 
& \makecell[l]{$\text{Step Size Ratio Error}_i$ \\ $(i = 1, \dots, N{-}1)$} & $\left| \frac{\text{step}_i}{\text{step}_{i+1}} - 2 \right| < \alpha \cdot \frac{1}{2^{N - i - 1} \cdot \sqrt{12N}}$ \\
& Sampling Error & $\left| V_{\text{ideal}} - V_{\text{sampled}} \right| < \alpha \cdot \frac{V_{\text{DD}}}{2^{N}\cdot\sqrt{12}}$ \\
& Thermal Noise & $\sqrt{\frac{2kT}{C} + \overline{v_{n,\text{cmp}}^2}} < \alpha \cdot \frac{V_{\text{DD}}}{2^{N}\cdot\sqrt{12}}$ \\
& Power & $\displaystyle min \ V_{\text{DD}} \cdot I_{\text{avg}}$ \\
\hline
\rule{0pt}{3.6ex}
\multirow{5}{*}{\makecell[c]{Accurate \\ Evaluation}} 
& SNDR & \makecell[l]{$10 \log_{10} \left( \frac{P_\text{signal}}{P_\text{tot, max}} \right)$ \\ $ < 6.02N - 4.25\ \text{dB}$} \\
& ENOB & $\frac{\text{SNDR} - 1.76}{6.02}$ \\
& FoM$_\text{W}$ & $\displaystyle min \ \frac{P}{2^{\text{ENOB}} \cdot f_\text{s}}$ \\
& FoM$_\text{S}$ & $\displaystyle min \ \text{SNDR} + 10 \log_{10} \left( \frac{f_\text{s}/2}{P} \right)$ \\
\hline
\end{tabular}}
\end{threeparttable}
\end{table}


\subsection{Low-Cost Simulation-Based Global Optimization}


Convergence can be slow without an effective optimization algorithm. Recent work has introduced machine learning into the optimization process. For example, \cite{budakEfficientAnalogCircuit2022} uses a Bayesian optimization framework in which an ANN-based surrogate model predicts the performance of offspring generated through differential evolution (DE). The use of ANN significantly reduces the needed simulations, but may introduce prediction errors. To mitigate this, only top candidates are simulated via infill sampling, with a Beta distribution-based ranking distinguishing good from poor designs. The method converges efficiently under stringent specifications for problems with fewer than 50 variables. This paper adopts the same components. The design variables in the optimization are:

\subsubsection{Bootstrap Switch}
The device sizes of the sampling switch are treated as design variables, as they directly affect both sampling accuracy and speed. To meet high-speed requirements, the sizes of the transistors in the clocking path are also parameterized, given their impact on switching speed.
\subsubsection{CDAC}
The unit capacitor size and CDAC driver transistor sizes are optimized. The driver size is scaled down progressively from the most significant bit to the least significant bit, typically halving at each step, until reaching the minimum size allowed by the PDK. This ensures proper settling behavior while minimizing power consumption.
\subsubsection{Comparator}
All its transistor sizes are included as design variables due to their critical role in determining noise characteristics, power consumption, and timing accuracy.
\subsubsection{DFF}
All transistor sizes are optimized, as they affect both timing precision and power efficiency.

To minimize simulation time, three tests are performed in parallel. The first measurement involves a single point transient simulation with the differential input fixed at $V_{DD}$. At the end of the sampling period, the sampling error can be directly evaluated. To check SSRE, voltage steps are measured during each bit cycle, and step ratios are recorded. The second and third measurements separately simulate the noise contributions from the CDAC and the comparator. Together, the RMS value is used to estimate the total thermal noise with high accuracy.


The global optimization completes once a specified number of parameters have converged. Convergence in the global phase is defined by a sufficiently small variance in the current population. Then the local optimizer starts from $\mathbf{x}_{\text{best}}$ and runs until full convergence. 

\subsection{Fast Local Optimization Using Parallel Multi-Fidelity Transient Simulation}
Local optimizers guarantee efficiency in finding a nearby optima in the vicinity of a good global starting point. However, they often face the 'curse of dimensionality', where the computational complexity grows exponentially as the number of design variables increases. To deal with this, in the local optimization stage, the converged design variables are fixed and only unconverged ones will be further optimized in this stage to prevent unnecessary move. The forms of local optimization algorithms vary depending on whether the problem is constrained and whether derivative information is available. To locally solve the constrained derivative-free optimization problem, the pattern search method \cite{moserHookeJeevesBasedMemetic2009} is adapted and described as in \textbf{Algorithm 1}. At each iteration, the algorithm explores predefined directions (e.g., along dimensions), accepting a new point if it improves the objective. Otherwise, it reduces the step size and continues until convergence. For this task, it typically requires 40-50 iterations for convergence.

\begin{algorithm}
\footnotesize
\caption{Blended Hooke-Jeeves with Selective Rollback}
\begin{algorithmic}[1]
\State \textbf{Init:} $x_0$, step sizes $\Delta_i$, frozen mask $M$, tolerance $\varepsilon$
\State $x_{\text{best}} \gets x_0$, $x_{\text{backup}} \gets x_0$, $w \gets 0.5$, $c \gets 0$
\State Evaluate $f_{\text{cheap}}(x_{\text{best}})$ and $f_{\text{expensive}}(x_{\text{best}})$
\State $f_{\text{backup}} \gets f_{\text{expensive}}(x_{\text{best}})$
\For{$k = 1$ to max\_iter}
    \State Exploratory search around $x_{\text{best}}$ (skip frozen)
    \If{improved $x_{\text{new}}$ found}
        \State Pattern move: extrapolate while improving
        \State $x_{\text{best}} \gets x_{\text{curr}}$, $c \gets c + 1$
        \If{$c \bmod \lambda = 0$}
            \State Evaluate $f_{\text{expensive}}(x_{\text{best}})$
            \State $\text{penalty} \gets a \cdot \max(0, f_{\text{exp}} - f_{\text{backup}})$
            \State $f_{\text{blend}} \gets (1 - w)f_{\text{cheap}} + w \cdot \text{penalty}$
            \If{$f_{\text{blend}} > f_{\text{cheap}}$}
                \State Rollback to $x_{\text{backup}}$, shrink $\Delta[\sim M]$
                \State $w \gets \min(w + \delta_w, 1)$
            \Else
                \State Update backup
            \EndIf
        \EndIf
    \Else
        \State Shrink $\Delta_i$ for non-improved, unfrozen $i$
    \EndIf
    \If{$\|\Delta[\sim M]\| < \varepsilon$} \textbf{break} \EndIf
\EndFor
\State \textbf{Return:} $x_{\text{best}}, f_{\text{cheap}}, f_{\text{expensive}}$
\label{algorithm}
\end{algorithmic}
\begin{tikzpicture}
    \draw[semithick] (0,0) -- (\linewidth,0);
\end{tikzpicture}
{\footnotesize

\textit{Notes:}
$\lambda$ is the frequency of expensive evaluations;
$a$ is the penalty scale factor;
$\delta_w$ is the update step for $w$ after rollback;
$\varepsilon$ is the termination tolerance;
$w$ is the weight in the blended cost function;
$\Delta_i$ are the coordinate step sizes;
$M$ is the frozen dimension mask.
}
\end{algorithm}

Accurate SNDR characterization of high-resolution SAR ADCs through transient simulation is often computationally intensive, due to the long simulation times required to capture sufficient output samples at the Nyquist rate. As the resolution increases, longer simulation durations and finer time steps are needed. These computational challenges not only slow down the verification process but also make it impractical to include conventional SNDR evaluation within design automation or optimization loops, where rapid and repeated performance assessments are essential.

\begin{figure}[tb]
\centering
\includegraphics[width=3in]{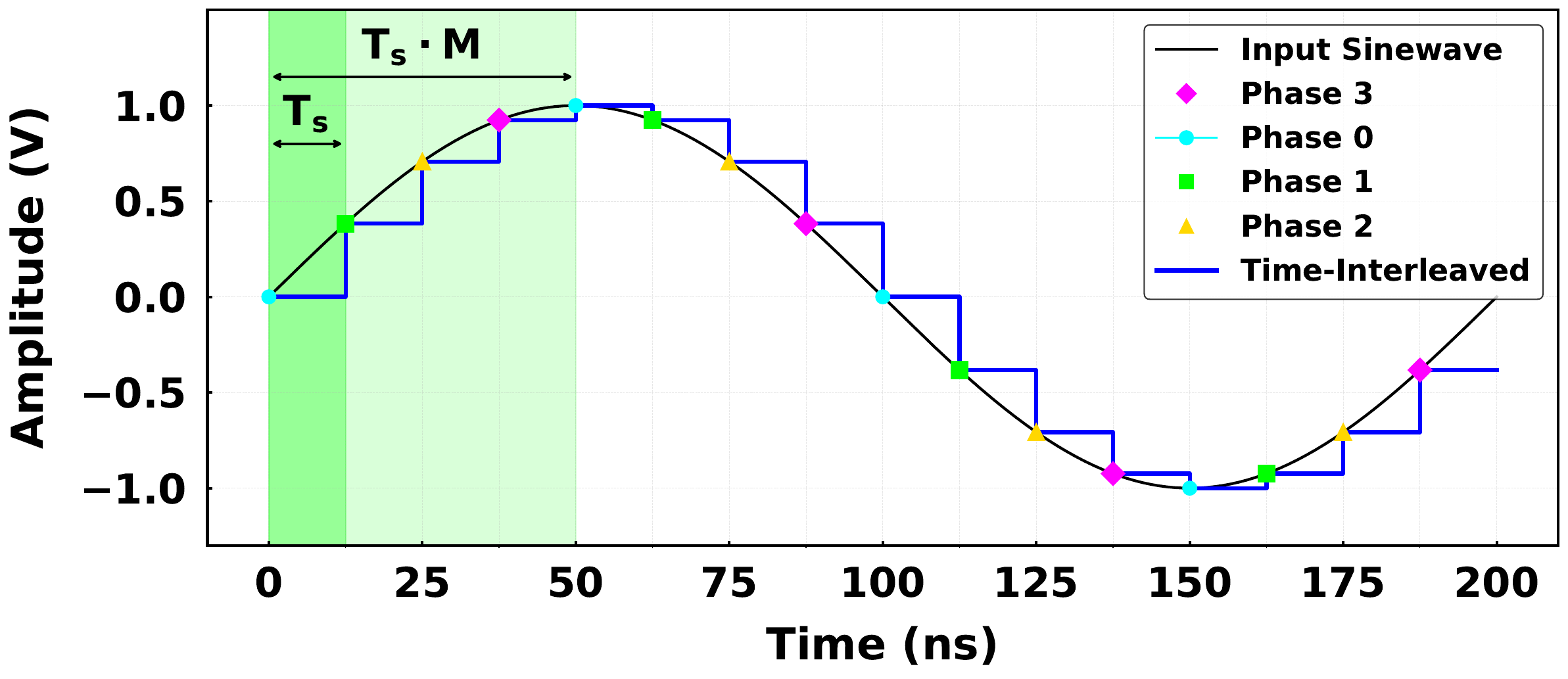}
\caption{Illustration of phase-shifted and time-interleaved parallel transient simulation: 16-point coverage via 4×4 samples.}
\label{sine}
\end{figure}

To mitigate this, we propose an efficient and scalable simulation method that leverages the periodic and deterministic nature of the input sinewave and the ADC’s response. Instead of performing a single long simulation at the full ADC sampling rate $f_s$, the method divides the task into $M$ parallel simulations, each operating at a reduced sampling rate of $f_s/M$. In each of the $M$ simulations, the input sine wave is phase-shifted by $\Delta \phi_k = \frac{2\pi k}{M}$, or equivalently, a time shift of $\Delta t_k = \frac{T_{\text{in}}}{M} \cdot k$, where $T_{\text{in}}$ is the sine wave period and $k = 0, 1, \ldots, M-1$. This ensures full phase coverage of the sine wave cycle across all simulations, while each simulation only needs to capture a reduced portion of the full-rate behavior. For instance, as shown in Fig. \ref{sine}, each parallel simulation collects 4 samples of the 4-bit ADC output data, which are then interleaved into a 16-point output data for FFT-related performance extraction. This method enables fast and accurate SNDR evaluation for high-resolution ADCs, specifically:
\begin{itemize}
    \item Reduced simulation time and memory usage without compromising SNDR accuracy.
    \item Scalable simulation speed-up proportional to $M$, via parallel execution.
    \item General applicability to any N-bit SAR ADC architecture.
    \item Easily integrable into design optimization and verification loops, where rapid and repeated FFT-based evaluations are required.
\end{itemize}

The proposed parallel simulation method is embedded within the local optimization algorithm. Even with a parallel setup, a single simulation can take several minutes. So in this framework, full-accuracy simulations are periodically invoked using the proposed technique to ensure accurate evaluation of SNDR. For intermediate query points during the local search, lower-cost coarse simulations in Section IV are applied. A blended objective function \( f_{\text{blended}}(x) \) is developed, which linearly combines \( f_{\text{cheap}}(x) \) and a penalty term derived from an expensive evaluation. The blending coefficient \( w \in [0, 1] \) determines the degree of trust in the penalty feedback: when \( w \) is small, the optimization follows the cheap evaluation, while a larger \( w \) shifts emphasis toward correcting misleading cheap evaluations. This adaptive mechanism balances efficiency with accuracy, ensuring both accurate convergence and practical runtime in the presence of local non-linearity.

\section{Experimental Results}
The sizing tool was developed using Matlab, with user inputs managed through a YAML file. The simulator is Cadence Spectre. The global optimization phase requires 3 hours on average, while the local optimization phase takes an additional 3 hours, using a 32-core machine running at 3.5 GHz.

\begin{figure}[t]
\centering
\subfloat[\label{3D}]{%
\includegraphics[width=0.5\linewidth]{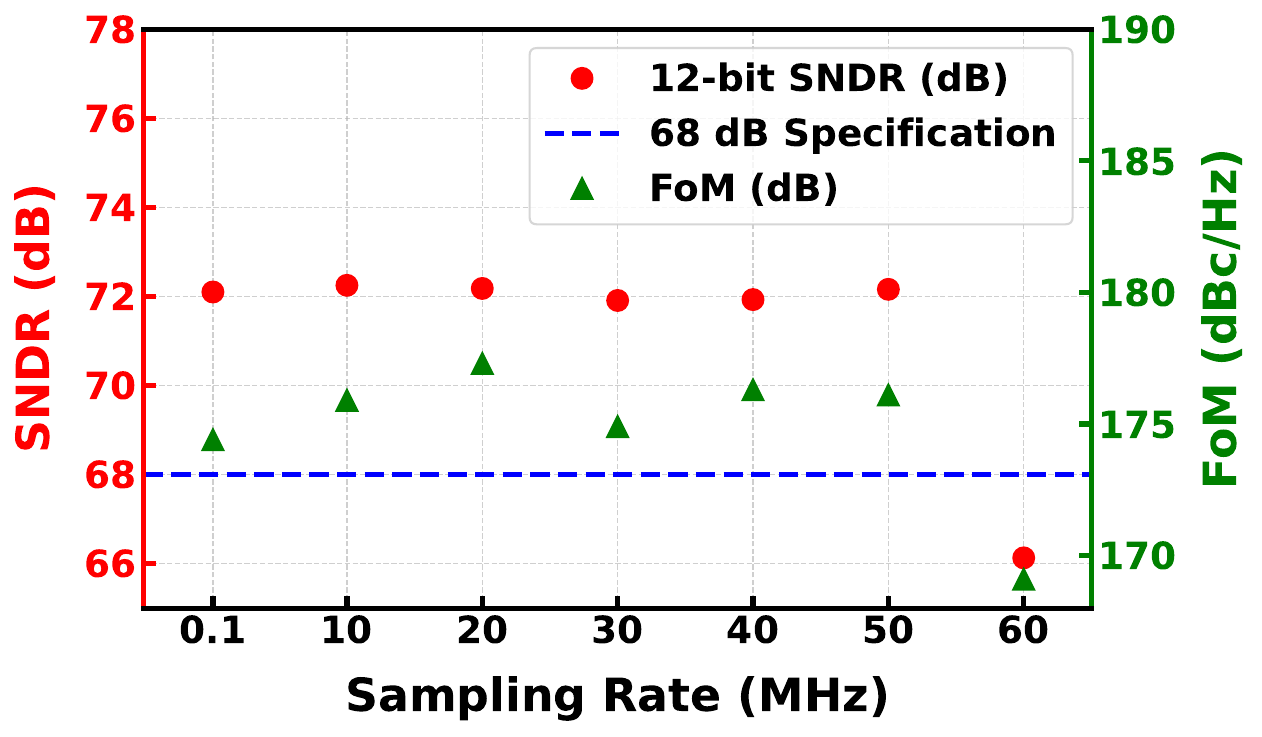}}
\hfill
\subfloat[\label{2D}]{%
\includegraphics[width=0.5\linewidth]{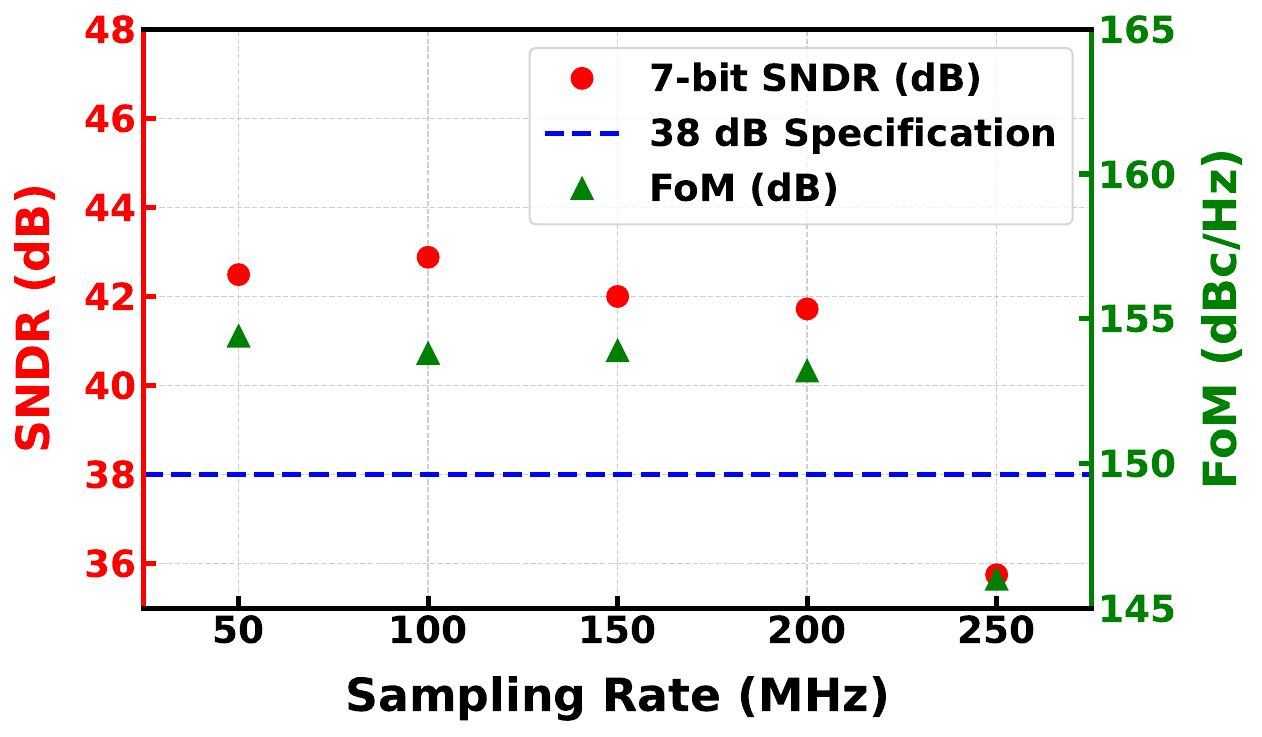}}
\caption{SNDR and FoM of 12 design cases: (a) 12 bit (b) 7 bit. 10 design cases with $\alpha$ = 1 and 2 with $\alpha$ = 2.}
\label{sndr}
\end{figure}

\begin{figure}[t]
\centering
\includegraphics[width=3.2in]{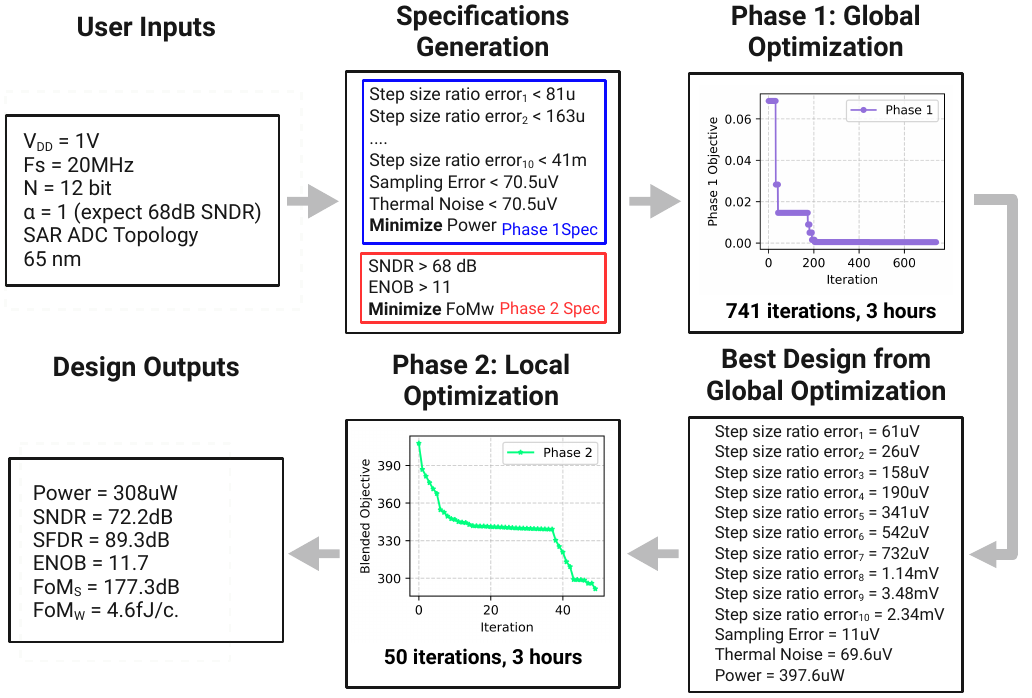}
\caption{Example sizing process for a 12 bit 20 MHz SAR ADC.}
\label{example}
\end{figure}

The proposed methodology was validated using a 65-nm CMOS process in 10 designs with $\alpha$ = 1 and 2 designs with $\alpha$ = 2, shown in Fig. \ref{sndr}. As illustrated, all design cases satisfy the specifications and achieve SNDR performance up to 72.2 dB and FoMs up to 177.3 dB. As an example, Fig. \ref{example} shows the optimization process of one design case.

\setlength{\tabcolsep}{2.2pt}  
\begin{table}[t]
\centering
\caption{Comparison with Prior SAR ADC Designs}
\label{comparison}
\renewcommand{\arraystretch}{1.1}
\begin{tabular}{l|cc|cc|cc|cc}
\toprule
\textbf{Feature} & 
\multicolumn{2}{c|}{\makecell{\textbf{TCAS-II-18\textsuperscript{†}} \\ \cite{seoReusableCodeBasedSAR2018}}} &
\multicolumn{2}{c|}{\makecell{\textbf{TVLSI-18\textsuperscript{†}} \\ \cite{dingHybridDesignAutomation2018}}} &
\multicolumn{2}{c|}{\makecell{\textbf{ICCAD-22\textsuperscript{‡}} \\ \cite{liuOpenSAROpenSource}}} &
\multicolumn{2}{c}{\textbf{This Work\textsuperscript{*}}} \\
\midrule
Process [nm] & 180 & 28 & 40 & 40 & 40 & 40 & \textbf{65} & \textbf{65} \\
Power Supply [V] & 1.8 & 1 & 1 & 1 & 1.2 & 0.7 & \textbf{1} & \textbf{1} \\
$F_s$ [MS/s] & 0.1 & 50 & 32 & 1 & 80 & 1 & \textbf{150} & \textbf{20} \\
Resolution [bit] & 12 & 11 & 8 & 12 & 10 & 12 & \textbf{7} & \textbf{12} \\
Power [$\mu$W] & 31.6 & 399 & 187 & 16.7 & 754.8 & 9.6 & \textbf{480} & \textbf{308} \\
SNDR [dB] & 63.3 & 56.8 & 47.4 & 61.1 & 56.3 & 68.8 & \textbf{42.0} & \textbf{72.2} \\
SFDR [dB] & 70.2 & 69.2 & 57.8 & 68.3 & 70.3 & 85.8 & \textbf{58.3} & \textbf{89.3} \\
ENOB & 10.2 & 9.1 & 7.6 & 9.9 & 9.1 & 11.1 & \textbf{6.68} & \textbf{11.7} \\
FoM$_\text{S}$ [dB]\textsuperscript{1} & 155.3 & 164.8 & 156.7 & 165.8 & 166.8 & 176.0 & \textbf{153.9} & \textbf{177.3} \\
FoM$_\text{W}$ [fJ/c.-step]\textsuperscript{2} & 265.5 & 14.1 & 30.7 & 18.1 & 10.8 & 4.3 & \textbf{31.6} & \textbf{4.6} \\
\bottomrule
\end{tabular}

\vspace{1mm}
\begin{minipage}{\linewidth}
\raggedright
\scriptsize
\textsuperscript{1} FoM$_\text{S}$ = SNDR + 10·log$_{10}$(Fs / 2 / Power).\\
\textsuperscript{2} FoM$_\text{W}$ = Power / (2$^{\text{ENOB}}$ · Fs).\\
\textsuperscript{†} Measurement results.\\
\textsuperscript{‡} Post-layout results.\\
\textsuperscript{*} Pre-layout results.
\end{minipage}
\vspace{-4mm}
\end{table}

Table \ref{comparison} summarizes and compares the proposed methodology and its performance with SAR ADC designs using block-level design approaches with similar resolutions and sampling rates. The proposed approach demonstrates superior speed performance compared to existing synthesized SAR ADC implementations. In addition, it significantly reduces the manual design effort by automating the whole design process, thereby enabling efficient exploration of the performance space. Although the methodology is demonstrated in 65 nm CMOS, the sizing framework is technology-independent, requiring only technology-specific inputs such as the PDK and design specifications. Furthermore, it is architecture-independent, as critical design decisions are made automatically by the tool rather than being manually specified by the designer.

\section{Conclusions\label{sec:conclusions}}
This brief presents a system-level global-local optimization framework for SAR ADC sizing, where the high-dimensional design space is first explored using a global optimizer and then refined through local optimization. By leveraging low-cost evaluations in the global phase and multi-fidelity simulations during local refinement, the proposed methodology efficiently manages overall sizing time. Experimental results demonstrate that the framework can achieve high-performance designs with minimal manual intervention across various resolutions and sampling rates.


\end{document}